\documentclass[graybox]{svmult}

\usepackage{mathptmx}       % selects Times Roman as basic font
\usepackage{helvet}         % selects Helvetica as sans-serif font
\usepackage{courier}        % selects Courier as typewriter font
\usepackage{type1cm}        % activate if the above 3 fonts are
\usepackage{makeidx}         % allows index generation
\usepackage{graphicx}        % standard LaTeX graphics tool
\usepackage{multicol}        % used for the two-column index
\usepackage[bottom]{footmisc}% places footnotes at page bottom

\usepackage{amsfonts}
\usepackage{amsmath, amssymb}
\usepackage{enumitem}

\newtheorem{assumption}{Assumption}

\newtheorem{Algorithm}{Algorithm}

\newcommand{\hp}{s}

\makeindex             % used for the subject index
\begin{document}

\title*{Diffuse interface approaches in atmosphere and ocean - modeling and numerical implementation}
\titlerunning{Diffuse interface approaches in atmosphere and ocean}

\author{Harald Garcke, Michael Hinze and Christian Kahle}
\institute{Harald Garcke \at Universit\"at Regensburg, Fakult\"at f\"ur Mathematik, 93040 Regensburg, \email{harald.garcke@mathematik.uni-regensburg.de} \and Michael Hinze \at Universit\"at Hamburg, Fachbereich Mathematik, Bundesstr. 55, 20146 Hamburg, \email{michael.hinze@uni-hamburg.de}
\and Christian Kahle \at Universit\"at Hamburg, Fachbereich Mathematik, Bundesstr. 55, 20146 Hamburg, \email{christian.kahle@uni-hamburg.de}}
\maketitle

\abstract{We propose to model physical effects at the sharp density interface
between atmosphere and ocean with the help of diffuse interface approaches
for multiphase flows with variable densities. We use the variable-density
model proposed in \cite{m6:AbelsGarckeGruen_CHNSmodell}.
This results in a Cahn-Hilliard/Navier-Stokes type system which we complement with tangential Dirichlet boundary conditions
to incorporate the effect of wind in the atmosphere.
Wind is responsible for waves at the surface of the ocean,
whose dynamics have an important impact on the $CO_2-$exchange between ocean and atmosphere.
We tackle this mathematical model numerically with fully adaptive and integrated
numerical schemes tailored to the simulation of variable density multiphase flows
governed by diffuse interface models. Here,
{\it fully adaptive, integrated, efficient, and reliable} means that the mesh resolution
is chosen by the numerical algorithm according to a prescribed error tolerance in the
{\it a posteriori} error control on the basis of residual-based error indicators,
which allow to estimate the true error from below (efficient) and from above (reliable).
Our approach is based on the work of
\cite{m6:HintermuellerHinzeKahle_adaptiveCHNS,m6:GarckeHinzeKahle_CHNS_AGG_linearStableTimeDisc},
where a fully adaptive
efficient and reliable  numerical method for the simulation of two-dimensional multiphase flows with variable densities is developed.
We incorporate the stimulation of surface waves via appropriate boundary conditions.
}

\section{Introduction}
\label{m6:sec:1}
The energy and momentum transfer from the atmosphere to the ocean is an essential ingredient for the accurate modelling of the energy cycle.  In fact, the vast majority of the energy input to the ocean comes from the winds ($\sim$20 TW), with much smaller inputs from tides ($\sim$3.5 TW) and geothermal heating ($\sim$0.05 TW) \cite[numbers taken from]{m6:wunsch2004}.  It is therefore not surprising that the energy transfer from the wind to the surface wave field, and the ensuing energy dissipation through breaking waves, represents the largest transfer of energy in the oceans \cite{m6:wunsch2004}.  Despite the enormous importance of the processes of surface wave generation and dissipation, there are still fundamental gaps in our ability to conduct both process modelling, and observational studies of these processes operating near the air-sea interface.

However, many major advances have recently been made through the use of powerful numerical simulations \cite[see the recent review by]{m6:sullivan2010}.  These simulations have shown that classical modelling and parameterisation techniques, such as the use of the law-of-the-wall turbulence scaling, must be revised to account for the dynamics of wind waves.  Here we outline a number of these results that must be included if an accurate, and energy consistent, treatment of atmosphere-ocean interactions is to be accomplished.

The transport of momentum from the atmosphere to the ocean in all but the finest scale modelling is largely done through specification of a bulk drag coefficient.  This coefficient attempts to capture the complex fluid flows around the air-water interface, as well as the dissipation processes therein.  This momentum (and energy) flux must then be partitioned between a number of different processes such as wind-wave growth, turbulent kinetic energy and dissipation through wave breaking, and the generation of wave-induced currents such as Langmuir circulations and Stokes drift.  All of these processes have been found to be important for the distribution of turbulent kinetic energy and dissipation in the surface mixed layer of the ocean \cite{m6:sullivan2007}.  Wave breaking has been found to lead to turbulence levels that are two orders of magnitude larger at the near-surface than the often used law-of-the-wall scaling predicts \cite{m6:drennan1996,m6:sutherland2015}.  Breaking is also responsible for generating large increases in the mean flow (Stokes drift) of the wave-affected near surface layer, which then provides a mechanism for the development of strong Langmuir circulations within the surface mixed layer \cite{m6:mcwilliams1997}.  These Langmuir circulations have been found to then redistribute the high turbulent kinetic energy throughout the surface mixed layer of the ocean \cite{m6:sullivan2007}.  Not only are the Langmuir circulations responsible for the redistribution of turbulent kinetic energy in the mixed layer, but this process is also found to contribute to the generation of internal waves at the base of the mixed layer that transport energy into the deeper ocean \cite{m6:polton2008}.

% summarise the results so far from DNS:
We propose  to use the thermodynamically consistent diffuse interface model proposed in
\cite{m6:AbelsGarckeGruen_CHNSmodell} to model the air-water interface between atmosphere and ocean.
This model will be extended to produce a series of direct numerical simulations of wind generated waves.
Only very few studies have examined the evolution of wind-waves using such a fundamental
approach \cite{m6:kihara2007,m6:shen2003,m6:sullivan2000,m6:sullivan2002,m6:tsai2007,m6:tsai2015,m6:lubin2015}.
However, these studies often do not involve a proper coupling of the water surface and the air flow above.
For example, the water surface is often replaced by another simpler boundary condition,
such as an impermeable sinusoidal wall \cite{m6:shen2003}, or an uncoupled propagating
water wave solution \cite{m6:sullivan2000,m6:sullivan2002}.
The recent study of \cite{m6:lubin2015} has shown how powerful a direct coupling of the air
and water layers is for predicting turbulent air-entraining structures in the breaking
of surface waves.  Another study that utilises a fully coupled treatment of the air
and water layers is that of \cite{m6:tsai2015}.  They show the important result that
turbulent water flows are generated even under the conditions of non-breaking surface waves.
We believe that the diffuse interface methods developed for the Cahn-Hilliard/Navier-Stokes
system will provide an improved method to deal with the current short comings of
simulating a direct coupling of the air-water interface. We recall here that one
reason is its flexibility in the numerical treatment of topology changes, which
might occur in e.g.~breaking waves, and another reason is the mass-conserving
property of the approach, see e.g.~\cite{m6:HintermuellerHinzeKahle_adaptiveCHNS}.
On the long run it is planed to use the method to simulate
wind-wave growth and be compare the numerical results to laboratory experiments using the PIV technique to resolve the airflow and water surface elevation.

% diffuse interface methods background:
The diffuse interface method of treating the air-water interface
using the Cahn-Hilliard/Navier-Stokes (CHNS) system forms
a new approach to the model studies described so far.
We note, however, that there are several contributions to numerical approaches
to the simulation of multiphase flows in the {\it sharp} interface formulation.
Here we refer, e.g., to the book of \cite{m6:GrossReusken_NumMethTwoPhaseFlow}, the work of
\cite{m6:GanesanTobiska09} as well as the works of \cite{m6:BGN14,m6:TurekWan07}. A benchmark for sharp interface
approaches to the numerical simulation of rising bubble dynamics is proposed
by \cite{m6:Hysing_Turek_quantitative_benchmark_computations_of_two_dimensional_bubble_dynamics},
which is accomplished with diffuse interface simulations
by \cite{m6:Aland_Voigt_bubble_benchmark}. A review of the development of phase-field
models and their numerical methods for multi-component fluid flows with interfacial
phenomena is given by \cite{m6:Kim12}. In the context of mechanical engineering
and meteorological applications phase-field models for two-phase flows are
often referred to as the {\it two-fluid formulation}, see
e.g. \cite{m6:MSSP10}, and \cite{m6:DE98}, compare also the related {\it volume-of-fluid schemes}, see e.g. \cite{m6:Lowengrub04}, as well as the references cited therein.

Since the dynamics of multiphase flows essentially depend on the dynamics
at the interfaces it is important to resolve the interfacial region
in diffuse interface models well. Here, adaptive numerical
concepts are the method of choice. Concerning the existing literature
on the solver development for the coupled CHNS system we note that
in \cite{m6:KayWelford_efficientsolution} a robust (with respect to the interfacial width)
nonlinear multigrid method was introduced with a double-well homogeneous
free energy density. We refer to \cite{m6:KayWelford_multigridsolver}
for the multigrid solver for the Cahn-Hilliard (CH) part only.
Later, in \cite{m6:KayStylesWelford} error estimates for the
coupled system were derived and numerically verified.
Coupled CHNS systems were also considered in \cite{m6:BoyerLapuerteMinjeaud}
with a double-well potential in the case of three-phase flows;
see also \cite{m6:Boyer_shear,m6:Boyer_two_phase_different_densities,m6:BoyerChupinFabrie}
and \cite{m6:AlandVoigt10,m6:Aland_Voigt_bubble_benchmark} as well as the references
therein for rather qualitative studies of the behaviour of multiphase and mixture flows.

Stable numerical schemes for  the recently developed thermodynamically consistent diffuse interface model
\cite{m6:AbelsGarckeGruen_CHNSmodell} are developed in
\cite{m6:GarckeHinzeKahle_CHNS_AGG_linearStableTimeDisc,m6:GruenMetzger__CHNS_decoupled,m6:Gruen_Klingbeil_CHNS_AGG_numeric}.
Concerning the numerical
treatment of the sole CH system many contributions can be found in the literature.
For a rather comprehensive discussion of available solvers we
refer to \cite{m6:HintermuellerHinzeTber}.
In the latter work, a {\it fully integrated} adaptive finite
element approach for the numerical treatment of the CH system
with a non-smooth homogeneous free energy density was developed.
The notion of a fully adaptive method relates to the fact that
the local mesh adaptation is based on rigorous residual
based \emph{a posteriori} error estimates (rather than heuristic
techniques based on, e.g., thresholding the discrete concentration or
the discrete concentration gradient).  The concept of an integrated
adaptation couples the adaptive cycle and the underlying solver
as the latter might need to be equipped with additional stabilisation
methodologies such as the Moreau-Yosida regularisation in the case
of non-smooth homogeneous free energy densities for guaranteeing mesh independence.
The latter is indeed obtained upon balancing regularisation and discretisation errors.
When equipped with a multi-grid scheme for solving the
linear systems occurring in the underlying semi-smooth Newton iteration,
an overall iterative scheme is obtained which is optimal in the sense
that the computational effort grows only linearly in the number of degrees of freedom.

%\bem{Own preliminary work in Diffuse Interface numerics approaches comes now...}
In \cite{m6:HintermuellerHinzeKahle_adaptiveCHNS} the approach of \cite{m6:HintermuellerHinzeTber}
is extended to a fully practical adaptive solver for the two-dimensional CHNS
system with a double obstacle potential according to \cite{m6:BloweyElliott_I}.
To the best of the applicants knowledge the work of \cite{m6:HintermuellerHinzeKahle_adaptiveCHNS}
contains the first rigorous approach to reliable and efficient residual based
\emph{a posteriori} error analysis for multi phase flows governed by diffuse
interface models. This approach is combined with a stable, energy conserving
time integration scheme in \cite{m6:GarckeHinzeKahle_CHNS_AGG_linearStableTimeDisc} to a fully reliable
and efficient adaptive and energy conserving \emph{a posteriori} concept
for the numerical treatment of variable density multiphase flows.
This approach was very successfully validated against the existing
sharp and diffuse interface rising benchmarks of
\cite{m6:Hysing_Turek_quantitative_benchmark_computations_of_two_dimensional_bubble_dynamics}
and \cite{m6:Aland_Voigt_bubble_benchmark}, respectively in the field.

\section{Diffuse interface approach}
\label{m6:sec:2}
\subsection{Notation}
Let $\Omega \subset \mathbb{R}^n$, $n\in \{2,3\}$ denote a bounded domain with boundary
$\partial \Omega$ and unit outer normal $\nu_\Omega$.
Let $I = (0,T]$ denote a time interval.

We use the conventional notation for Sobolev and Hilbert Spaces, see e.g.
\cite{m6:Adams_SobolevSpaces}.
With $L^p(\Omega)$, $1\leq p\leq \infty$, we denote the space of measurable functions on $\Omega$,
whose modulus to the power $p$ is Lebesgue-integrable. $L^\infty(\Omega)$ denotes the space of
measurable functions on $\Omega$, which are essentially bounded.
For $p=2$ we denote by
$L^2(\Omega)$ the space of square integrable functions on $\Omega$
with inner product $(\cdot,\cdot)$ and norm $\|\cdot \|$.
For a subset $D\subset \Omega$ and  functions $f,g\in L^2(\Omega)$ we by $(f,g)_D$  denote the inner product of
$f$ and $g$ restricted to $D$, and by $\|f\|_D$ the respective norm.
By $W^{k,p}(\Omega)$, $k\geq 1, 1\leq p\leq \infty$, we denote the Sobolev space of functions
admitting weak derivatives up to order $k$ in $L^p(\Omega)$.
If $p=2$ we write $H^k(\Omega)$.
The subset $H^1_0(\Omega)$ denotes $H^1(\Omega)$ functions
with vanishing boundary trace.\\
We further set
\begin{align*}
  L^2_{(0)}(\Omega) = \{v\in L^2(\Omega)\,|\, (v,1) = 0\},
\end{align*}
and with
\begin{align*}
  H(\mbox{div},\Omega) = \{v\in H^1_0(\Omega)^n\,|\, (\mbox{div}(v),q) = 0\,\forall q\in
  L^2_{(0)}(\Omega) \}
\end{align*}
we denote the space of all weakly solenoidal $H^1_0(\Omega)$ vector fields.\\
For $u\in L^q(\Omega)^n$, $q>n$, and $v,w \in H^1(\Omega)^n$ we introduce the
trilinear form
\begin{align*}
  a(u,v,w) =
  \frac{1}{2}\int_\Omega\left( \left(u\cdot\nabla\right)v \right)w\,dx
  -\frac{1}{2}\int_\Omega\left( \left(u\cdot\nabla\right)w \right)v\,dx.
\end{align*}
Note that there holds $a(u,v,w) = -a(u,w,v)$, and especially $a(u,v,v) = 0$.

\subsection{The mathematical model}
\label{m6:subsec:2.1}
In the present work we consider the following diffuse interface model for
two-phase flows with variable densities proposed in \cite{m6:AbelsGarckeGruen_CHNSmodell}:
\begin{align}
  \rho\partial_t v + \left( \left( \rho v + J\right)\cdot\nabla
  \right)v
  - \mbox{div}\left(2\eta Dv\right) + \nabla p =& \mu\nabla \varphi + \rho g &&
  \forall x\in\Omega,\, \forall t \in I,\label{m6:eq:CHNSstrong1}\\
  \mbox{div}(v) = &0&&
  \forall x\in\Omega,\, \forall t \in I,\label{m6:eq:CHNSstrong2}\\
  \partial_t \varphi + v \cdot\nabla \varphi - \mbox{div}(m\nabla \mu) = &0&&
  \forall x\in\Omega,\, \forall t \in I,\label{m6:eq:CHNSstrong3}\\
  -\sigma\epsilon \Delta \varphi + \frac{\sigma}{\epsilon}F'(\varphi) - \mu =&
  0&& \forall x\in\Omega,\, \forall t \in I,\label{m6:eq:CHNSstrong4}\\
  v(0,x) =& v_0(x)      &&    \forall x \in \Omega,\label{m6:eq:CHNSstrongIC1}\\
  \varphi(0,x) =& \varphi_0(x) &&\forall x \in \Omega,\label{m6:eq:CHNSstrongIC2}\\
  v(t,x) =& 0 &&\forall x \in \partial \Omega,\, \forall t \in I,\label{m6:eq:CHNSstrongBC1}\\
  \nabla \mu(t,x)\cdot \nu_\Omega =
  \nabla \varphi(t,x) \cdot \nu_\Omega =& 0  &&\forall x \in \partial
  \Omega,\, \forall t \in I,\label{m6:eq:CHNSstrongBC2}
\end{align}
where $J = -\frac{d\rho}{d\varphi}m\nabla \mu$.
Here $\Omega \subset \mathbb{R}^n,\, n\in \{2,3\}$, denotes an open and bounded domain, $I=(0,T]$
with $0<T<\infty$ a time interval,
$\varphi$ denotes the phase field,
$\mu$ the chemical potential,
$v$ the volume averaged velocity,
$p$ the pressure,
and $\rho = \rho(\varphi) = \frac{1}{2}\left((\rho_2-\rho_1)\varphi +(\rho_1+\rho_2)\right)$ the mean density,
where $0<\rho_1\leq\rho_2$ denote the densities of the involved fluids.
The viscosity is denoted by $\eta$ and can be chosen as an arbitrary positive function fulfilling $\eta(-1) =
\tilde\eta_1$ and $\eta(1) = \tilde \eta_2$, with individual fluid viscosities $\eta_1,\eta_2$.
The mobility is denoted by $m = m(\varphi)$. The gravitational force is denoted
by $g$.
By $Dv = \frac{1}{2}\left(\nabla v + (\nabla v)^t\right)$ we denote the symmetrized gradient.
The scaled surface tension is denoted by $\sigma$ and the interfacial
width is proportional to $\epsilon$.
The free energy is denoted by $F$.
For $F$ we use a splitting $F = F_+ + F_-$, where $F_+$ is convex and $F_-$
is concave.

The above model couples the Navier--Stokes equations \eqref{m6:eq:CHNSstrong1}--\eqref{m6:eq:CHNSstrong2}
to the Cahn--Hilliard model \eqref{m6:eq:CHNSstrong3}--\eqref{m6:eq:CHNSstrong4} in a thermodynamically
consistent way, i.e. a free energy inequality holds.
It is the main goal to introduce and analyze an (essentially) linear time discretization scheme for
the numerical treatment of  \eqref{m6:eq:CHNSstrong1}--\eqref{m6:eq:CHNSstrongBC2}, which also on the
discrete level fulfills the free energy inequality. This in conclusion leads to a stable scheme
that is thermodynamically consistent on the discrete level.

Existence of weak solutions to system \eqref{m6:eq:CHNSstrong1}--\eqref{m6:eq:CHNSstrongBC2} for a
specific class of free energies $F$ is shown in
\cite{m6:AbelsDepnerGarcke_CHNS_AGG_exSol,m6:AbelsDepnerGarcke_CHNS_AGG_exSol_degMob}.
See also the work \cite{m6:Gruen_convergence_stable_scheme_CHNS_AGG}, where the existence of
weak solutions for a different class of free energies $F$ is shown by passing to the limit in a
numerical scheme.
We refer to
\cite{m6:AkiDreyer_QuasiIncompressibleInterfaceModel},
\cite{m6:Boyer_two_phase_different_densities},
\cite{m6:Ding_Spelt_Shu_diffuse_interface_model_diff_density},
\cite{m6:Lowengrub_CahnHilliard_and_Topology_transitions},
and the review \cite{m6:AndersenFaddenWheeler} for other diffuse interface models for two-phase
incompressible flow.
Numerical approaches for different variants of the Navier--Stokes Cahn--Hilliard system have been
studied in
\cite{m6:Aland_Voigt_bubble_benchmark},
\cite{m6:Boyer_two_phase_different_densities},
\cite{m6:Feng_FullyDiscreteNSCH},
\cite{m6:Gruen_convergence_stable_scheme_CHNS_AGG},
\cite{m6:Gruen_Klingbeil_CHNS_AGG_numeric},
\cite{m6:GuoLinLowengrub_numericalMethodForCHNS_Lowengrub},
\cite{m6:HintermuellerHinzeKahle_adaptiveCHNS},
\cite{m6:HintermuellerHinzeKahle_adaptiveCHNS},
\cite{m6:GruenMetzger__CHNS_decoupled},
and
\cite{m6:KayStylesWelford}.

Our numerical treatment approach is based on the following weak formulation,
which is proposed in \cite{m6:GarckeHinzeKahle_CHNS_AGG_linearStableTimeDisc}.

\begin{definition}
We call $v$, $\varphi$, $\mu$ a
weak solution to \eqref{m6:eq:CHNSstrong1}--\eqref{m6:eq:CHNSstrongBC2}
if $v(0) = v_0$, $ \varphi(0) = \varphi_0$, $v(t) \in H(\mbox{div},\Omega)$ for $a.e.\,t \in I$ and
\begin{align}
  \frac{1}{2}\int_{\Omega}\left( \partial_t(\rho v)+\rho\partial_t v \right) w\,dx
  +\int_{\Omega}2\eta Dv:Dw\,dx&\nonumber\\
  +a(\rho v + J,v,w)
  =   \int_\Omega \mu\nabla \varphi w +\rho g w\,dx &\quad \forall w\in H(\mbox{div},\Omega),
  \label{m6:eq:CHNS1_weak}\\
  \int_\Omega\left(\partial_t\varphi + v\cdot \nabla \varphi\right) \Phi\,dx
  + \int_\Omega m(\varphi)\nabla \mu \cdot \nabla \Phi\,dx =0
  &\quad \forall \Phi \in  H^1(\Omega),
  \label{m6:eq:CHNS2_weak}\\
  \sigma \epsilon\int_\Omega\nabla \varphi\cdot\nabla \Psi\,dx
  +\frac{\sigma}{\epsilon}\int_\Omega F'(\varphi)\Psi\,dx
  - \int_\Omega \mu\Psi\,dx = 0 &\quad \forall \Psi \in H^1(\Omega),
  \label{m6:eq:CHNS3_weak}
\end{align}
is satisfied for almost all $t\in I$.
\end{definition}
For the assumptions on the data we refer to \cite{m6:GarckeHinzeKahle_CHNS_AGG_linearStableTimeDisc}. In the present work we
use the relaxed double-obstacle free energy given by
\begin{align}
  F(\varphi) & = \frac{1}{2}\left(1-\varphi^2+\hp\lambda^2( \varphi)\right),
  \label{m6:eq:F_rel_do}
\end{align}
with
\begin{align*}
  \lambda(\varphi)& := \max(0,\varphi-1) + \min(0,\varphi+1),
\end{align*}
where $\hp\gg 0$ denotes the relaxation parameter.
$F$  is introduced in \cite{m6:HintermuellerHinzeTber} as Moreau--Yosida
relaxation of the double-obstacle free energy
\begin{align*}
  F^{obst}(\varphi) =
  \begin{cases}
\frac{1}{2}\left(1-\varphi^2\right) & \mbox{ if } |\varphi|\leq 1,\\
0           & \mbox{ else},
\end{cases}
\end{align*}
which is proposed in \cite{m6:BloweyElliott_I} to model phase separation.

Let $v,\varphi,\mu$ be a sufficiently smooth
solution to
\eqref{m6:eq:CHNS1_weak}--\eqref{m6:eq:CHNS3_weak}.
Then we have from \cite{m6:GarckeHinzeKahle_CHNS_AGG_linearStableTimeDisc}  the energy relation
\begin{align}\label{m6:Energyrelation}
  \frac{d}{dt}
  \left(
  \int_\Omega  \frac{\rho}{2} |v|^2
  + \frac{\sigma\epsilon}{2}|\nabla \varphi|^2
  + \frac{\sigma}{\epsilon}F(\varphi)\,dx
  \right)
  = -\int_\Omega 2\eta |Dv|^2 + m|\nabla \mu|^2\,dx
  + \int_\Omega  \rho g  v\,dx .
\end{align}

\section{Discretization}
\label{m6:subsec:2.2}
\subsection{The time discrete setting}\label{m6:sec:timeDisc}

We now introduce a time discretization which mimics the energy inequality in \eqref{m6:Energyrelation} on the discrete level.
Let
$0=t_0<t_1<\ldots<t_{k-1}<t_k<t_{k+1}<\ldots<t_M=T$
denote an equidistant subdivision of the
interval $\overline I = [0,T]$ with $\tau_{k+1}-\tau_k = \tau$.
From here onwards the superscript $k$ denotes the corresponding variables at time instance
$t_k$.

\bigskip

\noindent \textbf{Time integration scheme}\\
Let $\varphi_0 \in H^1(\Omega)$
% with $|\varphi_0|\leq 1 \,a.e.$
and $v_0 \in  H(\mbox{div},\Omega)$.

\medskip

\noindent \textit{Initialization for $k=0$:}\\
Set $\varphi^0 = \varphi_0$ and $v^0=v_0$.\\
Find $\varphi^1 \in H^1(\Omega)$, $\mu^1\in H^1(\Omega)$, $v^1 \in H(\mbox{div},\Omega)$, such that
for all $w\in H(\mbox{div},\Omega)$, $\Phi \in H^1(\Omega)$, and $\Psi \in H^1(\Omega)$ it holds
\begin{align}
  \frac{1}{\tau}\int_\Omega \rho^1(v^1-v^0) w \,dx
  +\int_\Omega ((\rho^0v^0+J^1)\cdot \nabla) v^1 \cdot w\,dx&\nonumber\\
  +\int_{\Omega} 2\eta^1  Dv^{1}:Dw\,dx
  -\int_\Omega \mu^{1}\nabla \varphi^1 w + \rho^1 g w\,dx  &= 0&&
  \forall w \in H(\mbox{div},\Omega),\label{m6:eq:TD:chns1_solenoidal_init}\\
  \frac{1}{\tau}\int_\Omega (\varphi^{1}-\varphi^0) \Phi \,dx +
  \int_\Omega(v^{0}\cdot \nabla \varphi^{0}) \Phi \, dx& \nonumber \\
  +  \int_\Omega m(\varphi^0)\nabla \mu^{1}\cdot\nabla \Phi\,dx &=0 &&\forall
  \Phi\in H^1(\Omega),\label{m6:eq:TD:chns2_solenoidal_init}\\
  \sigma \epsilon\int_\Omega\nabla \varphi^{1}\cdot\nabla \Psi\,dx
  - \int_\Omega \mu^{1}\Psi\,dx \nonumber \\
  + \frac{\sigma}{\epsilon}\int_\Omega
  ((F_+)'(\varphi^{1})+(F_-)'(\varphi^0))\Psi\,dx
  &= 0&&\forall \Psi\in H^1(\Omega),\label{m6:eq:TD:chns3_solenoidal_init}
\end{align}
where $ J^1 := -\frac{d\rho}{d\varphi}(\varphi^1)m^1\nabla \mu^1$.

\medskip

\noindent \textit{Two-step scheme for $k\geq 1$:} \\
Given $\varphi^{k-1}\in H^1(\Omega)$,
$\varphi^k\in H^1(\Omega)$,
$\mu^k \in W^{1,q}(\Omega)$, $q>n$,
$v^k \in H(\mbox{div},\Omega)$,\\
find
$v^{k+1} \in H(\mbox{div},\Omega)$, $\varphi^{k+1}\in H^1(\Omega)$, $\mu^{k+1}\in
H^{1}(\Omega)$ satisfying
\begin{align}
  \frac{1}{2\tau}\int_\Omega \left( \rho^kv^{k+1}-\rho^{k-1}v^{k} \right) w +
  \rho^{k-1}(v^{k+1}-v^k)w\,dx&\nonumber\\
  +a(\rho^k v^k+J^k,v^{k+1},w)+\int_{\Omega} 2\eta^k Dv^{k+1}:Dw\,dx\nonumber\\
  -   \int_\Omega \mu^{k+1}\nabla \varphi^k w - \rho^{k} g w\,dx &= 0&&
  \forall w \in H(\mbox{div},\Omega),\label{m6:eq:TD:chns1_solenoidal}\\
  \frac{1}{\tau}\int_\Omega (\varphi^{k+1}-\varphi^k) \Phi \,dx +
  \int_\Omega(v^{k+1}\cdot \nabla \varphi^k) \Phi \, dx& \nonumber \\
  +  \int_\Omega m(\varphi^k)\nabla \mu^{k+1}\cdot\nabla \Phi\,dx &=0 &&\forall
  \Phi\in H^1(\Omega),\label{m6:eq:TD:chns2_solenoidal}\\
  \sigma \epsilon\int_\Omega\nabla \varphi^{k+1}\cdot\nabla \Psi\,dx
  - \int_\Omega \mu^{k+1}\Psi\,dx \nonumber \\
  + \frac{\sigma}{\epsilon}
  \int_\Omega  ((F_+)'(\varphi^{k+1})+(F_-)'(\varphi^k))\Psi\,dx
  &= 0&&\forall \Psi\in H^1(\Omega),\label{m6:eq:TD:chns3_solenoidal}
\end{align}
where $ J^k := -\frac{d\rho}{d\varphi}(\varphi^k)m^k\nabla \mu^k$.

We note that in \eqref{m6:eq:TD:chns1_solenoidal}--\eqref{m6:eq:TD:chns3_solenoidal} the only nonlinearity
arises from $F_+'$ and thus only the equation \eqref{m6:eq:TD:chns3_solenoidal} is nonlinear.
For a discussion of this scheme we refer to \cite{m6:GarckeHinzeKahle_CHNS_AGG_linearStableTimeDisc}. In
\cite{m6:Gruen_Klingbeil_CHNS_AGG_numeric} Gr\"un and Klingbeil propose a time-discrete solver
for \eqref{m6:eq:CHNSstrong1}--\eqref{m6:eq:CHNSstrongBC2} which leads to strongly coupled systems for
$v,\varphi$ and $p$ at every time step and requires a fully nonlinear solver. For this scheme Gr\"un
in \cite{m6:Gruen_convergence_stable_scheme_CHNS_AGG} proves an energy inequality and the existence of
so called generalized solutions.

\subsection{The fully discrete setting and energy inequalities}\label{m6:sec:fullyDiscrete}
For a numerical treatment we next discretize the weak formulation
\eqref{m6:eq:TD:chns1_solenoidal}--\eqref{m6:eq:TD:chns3_solenoidal} in space.
We aim at an adaptive discretization of the domain $\Omega$, and thus to have a different
spatial discretization in every time step.

Let $\mathcal T^{k} = \bigcup_{i=1}^{NT}T_i$ denote a conforming triangulation of
$\overline \Omega$ with closed simplices $T_i,i=1,\ldots,NT$ and edges $E_i,i=1,\ldots,NE$,
$\mathcal{E}^{k} = \bigcup_{i=1}^{NE}E_i$. Here $k$ refers to the time instance $t_{k}$.
On $\mathcal T^{k}$ we define the following finite element spaces:
\begin{align*}
  \mathcal{V}^{1}(\mathcal T^{k})
  =& \{v\in C(\mathcal{T}^{k}) \, |
  \, v|_T \in P^1(T)\, \forall T\in  \mathcal{T}^{k}\}
  =: \mbox{span}\{\Phi^i\}_{i=1}^{NP},\\
  \mathcal{V}^{2}(\mathcal T^{k}) =& \{v\in C(\mathcal{T}^{k}) \, |
  \, v|_T \in P^2(T)\, \forall T\in  \mathcal{T}^{k}\},
\end{align*}
where $P^l(S)$ denotes the space of polynomials up to order $l$ defined on $S$.

We introduce the discrete analogon to the space $H(\mbox{div},\Omega)$:
\begin{align*}
  H(\mbox{div},\mathcal{T}^{k}) &=
  \{ v\in \mathcal V^2(\mathcal T^{k})^n\,|\, (\mbox{div}v,q) =  0\,
  \forall q\in \mathcal V^1(\mathcal T^{k}) \cap L^2_{(0)}(\Omega),\,
  v|_{\partial \Omega}=0\}\\
  &:= \mbox{span}\{b^i\}_{i=1}^{NF},
\end{align*}

We further  introduce a $H^1$-stable projection operator
$\mathcal{P}^{k} : H^1(\Omega) \to \mathcal{V}^1(\mathcal T^{k})$
satisfying
\begin{align*}
  \|\mathcal{P}^{k}v\|_{L^p(\Omega)} \leq \|v\|_{L^p(\Omega)}
  \mbox{  and  }
  \|\nabla \mathcal{P}^{k}v \|_{L^r(\Omega)} \leq \|\nabla v\|_{L^r(\Omega)}
\end{align*}
for $v \in H^1(\Omega)$ with $r \in [1,2]$ and $p\in [1,6)$ if $n=3$, and $p\in [1,\infty)$ if
$n=2$.
Possible choices are the $H^1$-projection, the Cl\'ement operator
(\cite{m6:Clement_Interpolation}) or, by restricting the preimage to $C(\overline \Omega) \cap H^1(\Omega)$, the Lagrangian interpolation operator.

\bigskip

Using these spaces we state the discrete counterpart of
\eqref{m6:eq:TD:chns1_solenoidal}--\eqref{m6:eq:TD:chns3_solenoidal}:

Let $k\geq 1$, given
$\varphi^{k-1}\in \mathcal{V}^1(\mathcal{T}^{k-1})$,
$\varphi^{k}\in \mathcal{V}^1(\mathcal{T}^{k})$,
$\mu^{k}\in \mathcal{V}^1(\mathcal{T}^{k})$,
$v^{k}\in  H(\mbox{div},\mathcal T^{k})$,
find
$v^{k+1}_h \in  H(\mbox{div},\mathcal{T}^{k+1})$,
$\varphi^{k+1}_h\in \mathcal V^1(\mathcal  T^{k+1})$,
$\mu^{k+1}_h \in \mathcal V^1(\mathcal T^{k+1})$
such that for all
$w \in  H(\mbox{div},\mathcal T^{k+1})$,
$\Phi \in \mathcal V^1(\mathcal T^{k+1})$,
$\Psi \in  \mathcal V^1(\mathcal T^{k+1})$
there holds:
\begin{align}
  \frac{1}{2\tau}(\rho^{k}v^{k+1}_h-\rho^{k-1}v^k+\rho^{k-1}(v^{k+1}_h-v^k),w)
  + a(\rho^kv^k+J^k,v_h^{k+1},w)\nonumber\\
  +(2\eta^kDv^{k+1}_h,D w)-(\mu^{k+1}_h\nabla\varphi^{k}+\rho^k g,w)
  &= 0,\label{m6:eq:FD:chns1_solenoidal}\\
  \frac{1}{\tau}(\varphi^{k+1}_h-\mathcal{P}^{k+1}\varphi^k,\Phi)+(m(\varphi^k)\nabla
  \mu^{k+1}_h,\nabla \Phi) +(v^{k+1}_h\nabla \varphi^k, \Phi)
  &=0,\label{m6:eq:FD:chns2_solenoidal}\\
  \sigma\epsilon(\nabla \varphi^{k+1}_h,\nabla  \Psi)+
  \frac{\sigma}{\epsilon}(F_+^\prime(\varphi^{k+1}_h)+F^\prime_-(\mathcal{P}^{k+1}\varphi^k),\Psi)
  -(\mu^{k+1}_h,\Psi) &=0, \label{m6:eq:FD:chns3_solenoidal}
\end{align}
where
$\varphi^0 = P\varphi_0$ denotes the $L^2$ projection of $\varphi_0$ in
$\mathcal{V}^1(\mathcal T^0)$,
$v^0 = P^Lv_0$ denotes the Leray projection of $v_0$ in $H(\mbox{div},\mathcal T^0)$ (see
\cite{m6:ConstantinFoias_NS}), and $\varphi_h^1,\mu_h^1,v_h^1$ are obtained from the fully discrete variant of
\eqref{m6:eq:TD:chns1_solenoidal_init}--\eqref{m6:eq:TD:chns3_solenoidal_init}.

We have from \cite{m6:GarckeHinzeKahle_CHNS_AGG_linearStableTimeDisc} that the fully discrete system
\eqref{m6:eq:FD:chns1_solenoidal}--\eqref{m6:eq:FD:chns3_solenoidal} admits a unique solution,
where the analysis crucially depends on an energy inequality for the solution $(\varphi^{k+1}_h$, $\mu_h^{k+1}$, $v_h^{k+1})$ of
\eqref{m6:eq:FD:chns1_solenoidal}--\eqref{m6:eq:FD:chns3_solenoidal}. The energy inequality, which is also proven in \cite{m6:GarckeHinzeKahle_CHNS_AGG_linearStableTimeDisc}, is given as:

For $k\geq 1$:
\begin{align}
  \frac{1}{2}\int_\Omega \rho^k\left|v^{k+1}_h\right|^2\,dx
  + \frac{\sigma\epsilon}{2}\int_\Omega |\nabla \varphi^{k+1}_h|^2\,dx
  + \frac{\sigma}{\epsilon}\int_\Omega F(\varphi^{k+1}_h)\,dx\nonumber\\
  + \frac{1}{2}\int_\Omega \rho^{k-1}|v^{k+1}_h-v^k|^2\,dx
  + \frac{\sigma\epsilon}{2}
  \int_\Omega |\nabla \varphi^{k+1}_h-\nabla\mathcal{P}^{k+1} \varphi^k|^2\,dx\nonumber\\
  + \tau\int_\Omega 2\eta^k|Dv^{k+1}_h|^2\,dx
  +\tau\int_\Omega m^k|\nabla \mu^{k+1}_h|^2\,dx \nonumber\\
  \leq
  \frac{1}{2}\int_\Omega \rho^{k-1}\left|v^{k}\right|^2\,dx
  + \frac{\sigma\epsilon}{2}\int_\Omega |\nabla \mathcal{P}^{k+1}\varphi^{k}|^2\,dx
  + \frac{\sigma}{\epsilon}\int_\Omega F(\mathcal{P}^{k+1}\varphi^{k})\,dx
  +\tau \int_\Omega  \rho^k g v^{k+1}_h.\label{m6:eq:FD:energyInequality}
\end{align}

In \eqref{m6:eq:FD:energyInequality} the Ginzburg Landau energy of the current phase
field $\varphi^{k+1}$ is estimated against the Ginzburg Landau energy of the projection of the old
phase field $\mathcal P^{k+1}(\varphi^k)$. Since our aim is to obtain global in time inequalities
estimating the energy of the new phase field against the energy of the old phase field at each time
step we assume

\begin{assumption}\label{m6:ass:FP_leq_F}
Let $\varphi^k\in \mathcal{V}^1(\mathcal{T}^{k})$ denote the phase field at time instance $t_k$.
Let $\mathcal P^{k+1}\varphi^k\in \mathcal{V}^1(\mathcal{T}^{k+1})$ denote the projection of
$\varphi^k$ in $\mathcal{V}^1(\mathcal{T}^{k+1})$. We assume that there holds
\begin{align}
  \frac{\sigma}{\epsilon}F(\mathcal P^{k+1}\varphi^k)
  + \frac{1}{2}\sigma\epsilon |\nabla \mathcal P^{k+1}\varphi^{k}|^2
  \leq
  \frac{\sigma}{\epsilon}F(\varphi^k)
  +\frac{1}{2}\sigma\epsilon |\nabla \varphi^{k}|^2. \label{m6:eq:FP_leq_F}
\end{align}
\end{assumption}
This assumption  means, that the Ginzburg Landau energy is not increasing through projection. Thus
no energy is numerically produced.

Assumption \ref{m6:ass:FP_leq_F} is in general not fulfilled for arbitrary sequences
$(\mathcal{T}^{k+1})$ of triangulations.
To ensure \eqref{m6:eq:FP_leq_F} a post processing step can be added to the adaptive space meshing, see
Section \ref{m6:sec:Adaptivity}.

With this assumption we immediately get

\begin{theorem}
Assume that for every $k=0,1,\ldots$ Assumption \eqref{m6:ass:FP_leq_F} holds.
Then for every $1\leq k < l$ we have
\begin{align*}
  \frac{1}{2}(\rho^{k-1}_hv^k_h,v^k_h) +&  \frac{\sigma}{\epsilon}\int_\Omega
  F(\varphi^k_h)\,dx + \frac{1}{2}\sigma\epsilon (\nabla \varphi^{k}_h,\nabla \varphi^k_h)
  + \tau\sum_{m=k}^{l-1}(\rho^{m}g,v_h^{m+1})\\
  &\geq
  \frac{1}{2}(\rho^{l-1}v^{l}_h,v^{l}_h) +  \frac{\sigma}{\epsilon}\int_\Omega
  F(\varphi^{l}_h)\,dx + \frac{1}{2}\sigma\epsilon (\nabla \varphi^{l}_h,\nabla \varphi^{l}_h)\\
  & + \sum_{m=k}^{l-1} (\rho^{m-1}(v^{m+1}_h-v^m_h),(v^{m+1}_h-v^m_h))\\
  & +\tau\sum_{m=k}^{l-1} (2\eta^{m}Dv^{m+1}_h,Dv^{m+1}_h)\\
  & +\tau\sum_{m=k}^{l-1} (m(\varphi^m_h)\nabla \mu ^{m+1}_h,\nabla \mu^{m+1}_h)\\
  & +\frac{1}{2}\sigma\epsilon
  \sum_{m=k}^{l-1}( \nabla \varphi^{m+1}_h-\nabla \mathcal{P}^{m+1}\varphi^{m}_h,
  \nabla \varphi^{m+1}_h-\nabla \mathcal{P}^{m+1}\varphi^{m}_h).
\end{align*}
\end{theorem}

Now we have shown that there exists a unique solution to
\eqref{m6:eq:FD:chns1_solenoidal}--\eqref{m6:eq:FD:chns3_solenoidal}. The energy
inequality can be used to obtain uniform bounds on the fully discrete solution. This in turn can be used to obtain a
solution to the time discrete system \eqref{m6:eq:TD:chns1_solenoidal}--\eqref{m6:eq:TD:chns3_solenoidal} by a Galerkin
method, see \cite{m6:GarckeHinzeKahle_CHNS_AGG_linearStableTimeDisc} for the details of the proof.

\begin{theorem}\label{m6:thm:TD:exSol}
Let $v^{k}\in  H(\mbox{div},\Omega)$,  $\varphi^{k-1}\in
H^1(\Omega)$, $\varphi^k\in H^1(\Omega)$, and $\mu^k\in W^{1,q}(\Omega),q>n$ be given data.
Then there exists a unique weak solution to \eqref{m6:eq:TD:chns1_solenoidal}--\eqref{m6:eq:TD:chns3_solenoidal}.
Moreover,
$\varphi^{k+1}\in H^2(\Omega)$ and $\mu^{k+1}\in H^2(\Omega)$ holds. Moreover, the following energy inequality is valid.

\begin{align*}
  \frac{1}{2}\int_\Omega \rho^k\left|v^{k+1}\right|^2\,dx
  + \frac{\sigma\epsilon}{2}\int_\Omega |\nabla \varphi^{k+1}|^2\,dx
  + \frac{\sigma}{\epsilon}\int_\Omega F(\varphi^{k+1})\,dx\nonumber\\
  + \frac{1}{2}\int_\Omega \rho^{k-1}|v^{k+1}-v^k|^2\,dx
  + \frac{\sigma\epsilon}{2}\int_\Omega |\nabla \varphi^{k+1}-\nabla \varphi^k|^2\,dx\nonumber\\
  + \tau\int_\Omega 2\eta^k|Dv^{k+1}|^2\,dx
  +\tau\int_\Omega m^k|\nabla \mu^{k+1}|^2\,dx \nonumber\\
  \leq
  \frac{1}{2}\int_\Omega \rho^{k-1}\left|v^{k}\right|^2\,dx
  + \frac{\sigma\epsilon}{2}\int_\Omega |\nabla \varphi^{k}|^2\,dx
  + \frac{\sigma}{\epsilon}\int_\Omega F(\varphi^{k})\,dx
  + \int_\Omega  \rho^kg v^{k+1}dx.
\end{align*}
\end{theorem}

In our presentation $F$ denotes the relaxed double-obstacle free energy depending on the
relaxation parameter $\hp$. Let $(v_\hp,\varphi_\hp,\mu_\hp)_{\hp\in  \mathbb{R}}$ denote the
sequence of solutions of \eqref{m6:eq:TD:chns1_solenoidal}--\eqref{m6:eq:TD:chns3_solenoidal} for a
sequence $(\hp_l)_{l\in\mathbb{N}}$. Then we are able to argue convergence to solutions of a
limit system related to the double obstacle free energy $F^{obst}$.
More specifically, from the linearity of \eqref{m6:eq:TD:chns1_solenoidal} and
\cite[Prop. 4.2]{m6:HintermuellerHinzeTber} we conclude, that there exists a subsequence, still
denoted by $(v_\hp,\varphi_\hp,\mu_\hp)_{\hp\in  \mathbb{R}}$, such that
\begin{align*}
  (v_\hp,\varphi_\hp,\mu_\hp)_{\hp\in  \mathbb{R}}
  \to (v^*,\varphi^*,\mu^*)\quad \mbox{ in }  H^1(\Omega),
\end{align*}
where $(v^*,\varphi^*,\mu^*)$ denotes the solution of
\eqref{m6:eq:TD:chns1_solenoidal}--\eqref{m6:eq:TD:chns3_solenoidal},
where $F^{obst}$ is chosen as free energy. Especially
$|\varphi^*|\leq 1$ holds. Details are given in \cite{m6:GarckeHinzeKahle_CHNS_AGG_linearStableTimeDisc}.

\subsection{A-Posteriori Error Estimation}\label{m6:sec:Adaptivity}
For an efficient solution of \eqref{m6:eq:FD:chns1_solenoidal}--\eqref{m6:eq:FD:chns3_solenoidal} we next describe
an a-posteriori error estimator based mesh refinement scheme that is reliable and efficient up to terms
of higher order and errors introduced by the projection. We also describe how
Assumption \ref{m6:ass:FP_leq_F} on the evolution  of the free energy, given  in
\eqref{m6:eq:FD:energyInequality}, under projection  is fulfilled in the discrete setting.

Let us briefly comment on available adaptive concepts for the spatial discretization of
Cahn--Hilliard/Navier--Stokes systems. Heuristic approaches exploiting knowledge of the location
of the diffuse interface can be found in
\cite{m6:KayStylesWelford,m6:Aland_Voigt_bubble_benchmark,m6:Gruen_Klingbeil_CHNS_AGG_numeric}.
In \cite{m6:HintermuellerHinzeKahle_adaptiveCHNS} a fully adaptive, reliable and efficient, residual
based error estimator for the Cahn--Hilliard part in the Cahn--Hilliard/Navier--Stokes system is
proposed, which extends the results of \cite{m6:HintermuellerHinzeTber} for Cahn--Hilliard to
Cahn--Hilliard/Navier--Stokes systems with Moreau--Yosida relaxation of the double-obstacle free
energy.
A residual based error estimator for Cahn--Hilliard systems with double-obstacle free energy is
proposed in  \cite{m6:BanasNurnbergAPosteriori}.

In the present section we propose a fully integrated adaptive concept for the fully coupled
Cahn--Hilliard/Navier--Stokes system, where we exploit the energy inequality of \eqref{m6:Energyrelation}.

For the numerical realization we switch to the primitive setting for the flow part of our equation system. The corresponding fully discrete system now reads:\\
For $k\geq 1$, given
$\varphi^{k-1}\in H^1(\Omega)$,
$\varphi^k\in H^1(\Omega)$,
$\mu^k \in W^{1,q}(\Omega),q>n$,
$v^k \in H^1_0(\Omega)^n$
find
$v^{k+1}_h \in  \mathcal{V}^2(\mathcal{T}^{k+1})$,
$p_h^{k+1}\in \mathcal V^1(\mathcal T^{k+1}),\, \int_\Omega p_h^{k+1}\,dx= 0$,
$\varphi^{k+1}_h\in \mathcal V^1(\mathcal  T^{k+1})$,
$\mu^{k+1}_h \in \mathcal V^1(\mathcal T^{k+1})$
such that for all
$w \in \mathcal V^2(\mathcal T^{k+1})$,
$q \in \mathcal V^1(\mathcal T^{k+1})$,
$\Phi \in \mathcal V^1(\mathcal T^{k+1})$,
$\Psi \in \mathcal V^1(\mathcal T^{k+1})$
there holds:
\begin{align}
  \frac{1}{2\tau}(\rho^{k}v^{k+1}_h-\rho^{k-1}v^k+\rho^{k-1}(v^{k+1}_h-v^k),w)
  + a(\rho^kv^k+J^k,v_h^{k+1},w)\nonumber\\
  +(2\eta^kDv^{k+1}_h,\nabla w)-(\mu^{k+1}_h\nabla\varphi^{k} + \rho^k g,w) - (p_h^{k+1},\mbox{div}
  w) &= 0,\label{m6:eq:AD:chns1_fullDisc}\\
  - (\mbox{div} v_h^{k+1},q) &= 0,\label{m6:eq:AD:chns2_fullDisc}\\
  \frac{1}{\tau}(\varphi^{k+1}_h-\mathcal{P}^{k+1}\varphi^k,\Phi)+(m(\varphi^k)\nabla
  \mu^{k+1}_h,\nabla \Phi) -(v^{k+1}_h \varphi^k,\nabla \Phi)
  &=0,\label{m6:eq:AD:chns3_fullDisc}\\
  \sigma\epsilon(\nabla \varphi^{k+1}_h,\nabla
  \Psi)+
  \frac{\sigma}{\epsilon}(F_+^\prime(\varphi^{k+1}_h)+F^\prime_-(\mathcal{P}^{k+1}\varphi^k),\Psi)-(\mu^{k+1}_h,\Psi) &=0. \label{m6:eq:AD:chns4_fullDisc}
\end{align}

Thus we use the famous Taylor--Hood LBB-stable $P2-P1$ finite
element for the discretization of the velocity - pressure field and piecewise linear and continuous
finite elements for the discretization of the phase field and the chemical potential. For other kinds of possible
discretizations of the velocity-pressure field we refer to e.g.
\cite{m6:Verfuerth_aPosteriori_new}.

Note that we perform integration by parts in \eqref{m6:eq:AD:chns3_fullDisc} in the
transport term, using the no-slip boundary condition for $v_h^{k+1}$. As soon as
$\mathcal{P}^{k+1}$ is a mass conserving projection we by testing equation \eqref{m6:eq:AD:chns3_fullDisc} with $\Phi = 1$ obtain the conservation of mass in the fully discrete scheme.

The link between equations
\eqref{m6:eq:AD:chns1_fullDisc}--\eqref{m6:eq:AD:chns4_fullDisc} and
\eqref{m6:eq:FD:chns1_solenoidal}--\eqref{m6:eq:FD:chns3_solenoidal} is established by the fact that for $v^{k+1}_h,\varphi^{k+1}_h,\mu^{k+1}_h$ denoting the unique  solution to
\eqref{m6:eq:FD:chns1_solenoidal}--\eqref{m6:eq:FD:chns3_solenoidal}, there exists a unique pressure
$p_h^{k+1} \in \mathcal{V}^1(\mathcal{T}^{k+1}),\, \int_\Omega p_h^{k+1}\,dx= 0$ such that
$(v^{k+1}_h,p_h^{k+1},\varphi_h^{k+1},\mu_h^{k+1})$ is a solution to \eqref{m6:eq:AD:chns1_fullDisc}--\eqref{m6:eq:AD:chns4_fullDisc}.
The opposite direction is obvious.

Next we describe the error estimator which we use in our computations. We follow \cite{m6:HintermuellerHinzeTber} and restrict the
presentation of its construction to the main steps.

We define the following error terms:
\begin{align*}
  e_v := &v^{k+1}_h - v^{k+1},&
  e_p := &p^{k+1}_h - p^{k+1},\\
  e_\varphi := &\varphi^{k+1}_h - \varphi^{k+1},&
  e_\mu := &\mu^{k+1}_h - \mu^{k+1},
\end{align*}
as well as the discrete element residuals
\begin{align*}
  r^{(1)}_h :=&
  \frac{\rho^{k}+\rho^{k-1}}{2}v^{k+1}_h   - \rho^{k-1}v^k
  + \tau (b^k\nabla)v^{k+1}_h
  +\frac{1}{2}\tau \mbox{div}(b^k)v^{k+1}_h \\
  &-2\tau \mbox{div}\left(\eta^kDv^{k+1}_h\right)
  + \tau \nabla p^{k+1}_h
  - \tau \mu^{k+1}_h\nabla \varphi^k -  \rho^k g,\\
  r^{(2)}_h := & \varphi^{k+1}_h-\mathcal{P}^{k+1}\varphi^k + \tau v^{k+1}_h\nabla \varphi^k
  - \tau \mbox{div}(m^k\nabla \mu^{k+1}_h),\\
  r^{(3)}_h  := & \frac{\sigma}{\epsilon} F'_+(\varphi^{k+1}_h) +
  \frac{\sigma}{\epsilon}F'_-(\mathcal{P}^{k+1}\varphi^k) - \mu^{k+1}_h,
\end{align*}
where $b^k := \rho^kv^k+J^k $.
Furthermore we define the error indicators
\begin{equation}
  \begin{aligned}
    \eta_T^{(1)} :=& h_T\|r^{(1)}_h\|_T, &
    \eta_E^{(1)} :=& h_E^{1/2}\|2\eta^k \left[Dv^{k+1}_h\right]_\nu \|_E,\\
    \eta_T^{(2)} := & h_T\|r^{(2)}_h\|_T, &
    \eta_E^{(2)} := & h_E^{1/2}\|m^k\left[\nabla \mu^{k+1}_h\right]_\nu\|_E, \\
    \eta_T^{(3)} := & h_T\|r^{(3)}_h\|_T, &
    \eta_E^{(3)} := & h_E^{1/2}\|\left[\nabla \varphi^{k+1}_h\right]_\nu\|_E.
  \end{aligned}
  \label{m6:eq:AD:etaT123_etaTE123}
\end{equation}
Here $\left[\cdot\right]_\nu$  denotes the jump of a discontinuous function in normal direction
$\nu$ pointing from the triangle with lower global number to the triangle with higher global
number. Thus $\eta_E^{(j)}$, $j=1,2,3$ measures  the jump of the corresponding variable
across the edge $E$, while $\eta_T^{(j)}$, $j=1,2,3$ measures the triangle wise residuals.

In \cite[Theorem 9]{m6:GarckeHinzeKahle_CHNS_AGG_linearStableTimeDisc} the following Theorem is proven.
\begin{theorem}\label{m6:thm:AD:errorEstimate}
There exists a constant $C>0$ only depending on the domain $\Omega$ and the regularity of the mesh
$\mathcal{T}^{k+1}$ such that
\begin{align*}
  \underline{\rho}\|e_v\|^2
  + \tau \underline\eta \|\nabla e_v\|^2
  +\tau \underline m\|\nabla e_\mu\|^2
  + \sigma\epsilon\|\nabla e_\varphi\|^2
  +
  \frac{\sigma}{\epsilon}(F'_+(\varphi^{k+1}_h)-F'_+(\varphi^{k+1}),e_\varphi)\\
  \leq C\left( \eta_\Omega^2 + \eta_{h.o.t} + \eta_C \right),
\end{align*}
holds with
\begin{align*}
  \eta_\Omega^2 = &
  \frac{1}{\tau\underline \eta}
  \sum_{T\in \mathcal{T}^{k+1}}
  \left(\eta_T^{(1)}\right)^2
  +
  \frac{\tau}{\underline\eta}\sum_{E\in \mathcal{E}^{k+1}}
  \left(\eta_E^{(1)}\right)^2\\
  &
  \frac{1}{\tau\underline m}
  \sum_{T\in \mathcal{T}^{k+1}}
  \left(\eta_T^{(2)}\right)^2
  +
  \frac{\tau}{\underline m}\sum_{E\in \mathcal{E}^{k+1}}
  \left(\eta_E^{(2)}\right)^2\\
  &
  \frac{1}{\sigma\epsilon}
  \sum_{T\in \mathcal{T}^{k+1}}
  \left(\eta_T^{(3)}\right)^2
  +
  \sigma\epsilon\sum_{E\in \mathcal{E}^{k+1}}
  \left(\eta_E^{(3)}\right)^2,\\
  \eta_{h.o.t.} =& \tau (\mbox{div}(e_v),e_p), \\
  \mbox{ and }\eta_C = &(\mathcal{P}^{k+1}\varphi^k-\varphi^k,e_\mu)
  -  \frac{\sigma}{\epsilon}(F'_-(\mathcal{P}^{k+1}\varphi^k) -
  F'_-(\varphi^k),e_\varphi).
\end{align*}
\end{theorem}

In the numerical part, this error estimator is used together with the mesh adaptation cycle
described in \cite{m6:HintermuellerHinzeTber}.
The overall adaptation cycle
\begin{center}
  SOLVE $\to$ ESTIMATE $\to$ MARK $\to$ ADAPT
\end{center}
is performed once per time step.
For convenience of the reader we state the marking strategy here.
\begin{Algorithm}[Marking strategy]
~\label{m6:alg:AD:marking}
\begin{itemize}
  \item Fix $a_{\min}>0$ and $a_{\max}>0$, and set $\mathcal{A} = \{T\in \mathcal{T}^{k+1}\,|\,
  a_{\min}\leq |T|\leq a_{\max}\}$.
  \item Define indicators:
  \begin{enumerate}
    \item $\eta_T = \frac{1}{\tau \underline \eta}\left(\eta_T^{(1)}\right)^2
    + \frac{1}{\tau \underline m}\left(\eta_T^{(2)}\right)^2
    + \frac{1}{\sigma\epsilon}\left(\eta_T^{(3)}\right)^2$,
    \item $\eta_{TE} =\sum_{E\subset T}\left[
    \frac{\tau}{\underline \eta}\left(\eta_{TE}^{(1)}\right)^2
    +\frac{\tau}{\underline m}\left(\eta_{TE}^{(2)}\right)^2
    +\sigma\epsilon \left(\eta_{TE}^{(3)}\right)^2\right]$.
  \end{enumerate}
  \item Refinement: Choose $\theta^r\in (0,1)$,
  \begin{enumerate}
    \item Find a set $R^T\subset\mathcal{T}^{k+1}$ with
    $\theta^r\sum_{T\in \mathcal{T}^{k+1}}\eta_T\leq \sum_{T \in R^T}\eta_T$,
    \item Find a set $R^{TE}\subset\mathcal{T}^{k+1}$ with
    $\theta^r\sum_{T\in \mathcal{T}^{k+1}}\eta_{TE}\leq \sum_{T \in R^{TE}}\eta_{TE}$.
  \end{enumerate}
  \item Coarsening: Choose $\theta^c\in (0,1)$,
  \begin{enumerate}
    \item Find the set $C^T\subset\mathcal{T}^{k+1}$ with
    $\eta_T\leq \frac{\theta^c}{N}\sum_{T \in \mathcal{T}^{k+1}}\eta_T\,\forall T\in C^T$,
    \item Find the set $C^{TE}\subset\mathcal{T}^{k+1}$ with
    $\eta_{TE}\leq \frac{\theta^c}{N}\sum_{T \in \mathcal{T}^{k+1}}\eta_{TE}\,\forall T\in C^{TE}$.
  \end{enumerate}
  \item Mark all triangles of $\mathcal{A}\cap (R^T\cup R^{TE})$ for refining.
  \item Mark all triangles of $\mathcal{A}\cap (C^T\cup C^{TE})$ for coarsening.
\end{itemize}
\end{Algorithm}

\subsubsection*{Ensuring the validity of the energy estimate}
To ensure the validity of the energy estimate during the numerical computations we  ensure
that Assumption \ref{m6:ass:FP_leq_F} holds triangle-wise.
For the following considerations we restrict to bisection as refinement strategy combined with the
$i$FEM coarsening strategy proposed in \cite{m6:iFEMpaper}. This strategy only coarsens patches
consisting of four triangles by replacing them by two triangles if the central node of the patch is
an inner node of $\mathcal{T}^{k+1}$, and patches consisting of two triangles by replacing them by
one triangle if the central node of the patch lies on the boundary of $\Omega$. A patch fulfilling
one of these two conditions we call a nodeStar.
By using this strategy,  we do not harm the Assumption \ref{m6:ass:FP_leq_F} on triangles that are
refined.
We note that this assumption can
only be violated on patches of triangles where coarsening appears.

After marking triangles for refinement and coarsening  and before applying refinement and
coarsening to $\mathcal{T}^{k+1}$ we make a post-processing of all triangles
that are marked for coarsening.

Let $M^C$ denote the set of triangles marked for coarsening obtained by the marking strategy
described in Algorithm \ref{m6:alg:AD:marking}.
To ensure the validity of the energy estimate \eqref{m6:eq:FD:energyInequality} we perform the
following post processing steps:
\begin{Algorithm}[Post processing]
~\label{m6:alg:AD:PostProcessing}
\begin{enumerate}
  \item For each triangle $T \in M^C$:\\
  if $T$ is not part of a nodeStar\\
  then set $M^C := M^C\setminus T$.
  \item For each nodeStar $S \in M^C$:\\
  if Assumption \eqref{m6:ass:FP_leq_F} is not fulfilled on $S$\\
  then set $M^C := M^C \setminus S$.
\end{enumerate}

\end{Algorithm}

The resulting set $M^C$ does only contain triangles yielding nodeStars on which the Assumption
\ref{m6:ass:FP_leq_F} is fulfilled.

\section{Numerics}\label{m6:sec:Numerics}
%section identifier: ..:num:..
Let us finally give a numerical example to show the applicability of
the provided method  to
the simulation of the complex interaction at the ocean - atmosphere interface.
We use the implementation  from
\cite{m6:GarckeHinzeKahle_CHNS_AGG_linearStableTimeDisc}, that was developed for
the validation of results for the rising bubble benchmark from
\cite{m6:Aland_Voigt_bubble_benchmark,m6:Hysing_Turek_quantitative_benchmark_computations_of_two_dimensional_bubble_dynamics}.
For this reason it is not adapted to the present situation and we will comment on the restrictions and the
future work to tackle the given problem after showing the numerical results.

\begin{figure}
  \centering
  \fbox{
  \includegraphics[width=0.8\textwidth]{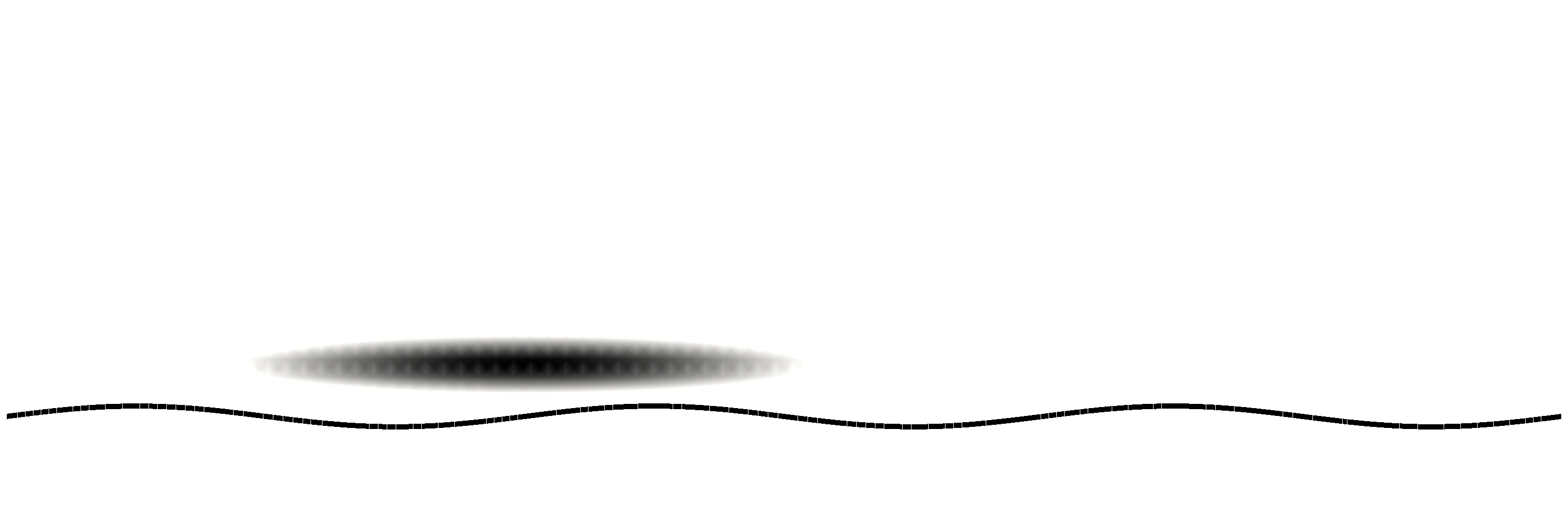}
  }
  \caption{The volume force $f$ that generates the 'wind' together with 
  the zero level line of $\varphi_0$.}
  \label{m6:fig:num:volForce}
\end{figure}

We use $\Omega = (0.0,3.0) \times (0.0,1.0)$
and a time horizon of $I = (0,10.0)$, that we subdivide into equidistant time
steps of length $\tau = 5e-4$. To mimic the wind forcing we introduce a volume
force $f = (f_1,f_2)^\top = (f_1,0)^\top$  as shown in Fig.
\ref{m6:fig:num:volForce} on the right hand side of the Navier--Stokes equation
 that is defined as
\begin{align*}
  (f_1(x)) =
  \begin{cases}
0 & \mbox{if } \|(x - m)/\sigma\| \geq 1,\\
\cos(\pi \|(x - m)/\sigma\|)^2 & \mbox{else},
\end{cases}
\end{align*}
where $m = (1.0,1.2)^\top$, $\sigma = (1.0,0.1)^\top$ and the division
$(x - m)/\sigma$ has to be understood component-wise. This is an approximation of a
gaussian bell with compact support.
In Fig.~\ref{m6:fig:num:volForce} we further show the zero level line of the
initial phase field $\varphi_0$, that is given by
\begin{align*}
  z(x_1,x_2):= & (x_2-0.02\sin(2\pi x_1) + 0.2)/\epsilon,\\
  z_0  = & \arctan{\sqrt{\hp-1}},\\
  \varphi_0(z) := &
  \begin{cases}
\sqrt{\frac{\hp}{\hp-1}} \sin(z)  & \mbox{ if } |z| \leq z_0,\\
\phantom{-}
\frac{1}{\hp-1}\left(s-\exp\left(\sqrt{s-1}(z_0-z)\right) \right)
& \mbox{ if }  z>z_0,\\
-\frac{1}{\hp-1}\left(s-\exp\left(\sqrt{s-1}(z_0+z)\right)\right) 
&  \mbox{if }z<-z_0.
\end{cases}
\end{align*}
Note that $z$ measures the distance in $x_2$ direction to the wave $0.2 -
0.02\sin(2\pi x_1)$ scaled by $\epsilon^{-1}$ and $\varphi_0(z)$ is the first order approximation to
a phase field with relaxed double-obstacle free energy \eqref{m6:eq:F_rel_do}, see
\cite[Sec. 10]{m6:kahle_dissertation}. The initial velocity is $v_0\equiv 0$ and
we have no-slip boundary data at $\partial \Omega$.

As parameters we choose $\rho_{water} = 1$, $\eta_{water} = 0.01$,
$\sigma_{water} = 0.0005$. Using unit velocity $V=1$ and unit length $d=1$, this
results in a Weber number of $We = 2000$ and after a required scaling due to
the chosen free energy, see \cite{m6:AbelsGarckeGruen_CHNSmodell},
we have $\sigma = 0.00032$. The gravity is $g = (0.0,-9.81)^\top$ and the mobility is $b = \frac{\epsilon}{500\sigma}$. We note that especially the chosen density $\rho_{water}$ does not correspond to the real world parameter. We use the air-water ratio
$\rho_{air} = 0.01\rho_{water}$ and $\eta_{air} = 0.01\eta_{water}$, which is ten times larger than the real world ratio. 
%Also these values do not correspond to physical values. 
To overcome the limitations of the
current implementation with respect to the density ratio is subject of future research.

For the adaptation process from Algorithm \ref{m6:alg:AD:marking} we choose
$\theta^r = 0.5$, $\theta^c = 0.01$, $V_{\min} = 8e-6$, $V_{\max} = 3e-4$.
For the numerical results presented we switched off the postprocessing proposed in Algorithm \ref{m6:alg:AD:PostProcessing}.
In \cite{m6:GarckeHinzeKahle_CHNS_AGG_linearStableTimeDisc} the influence of this postprocessing on the numerical simulation of the rising bubble benchmark is investigated in detail.

\begin{figure}
  \centering
  \hfill
  \includegraphics[width = 0.49\textwidth]{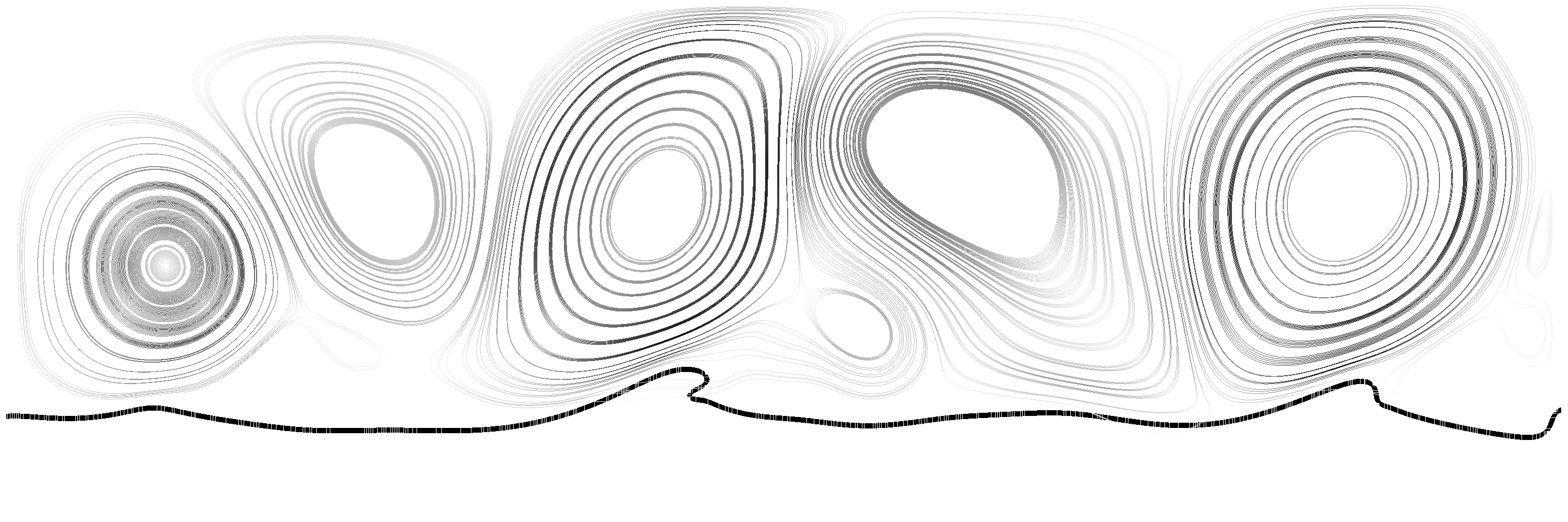}
  \hfill
  \includegraphics[width = 0.49\textwidth]{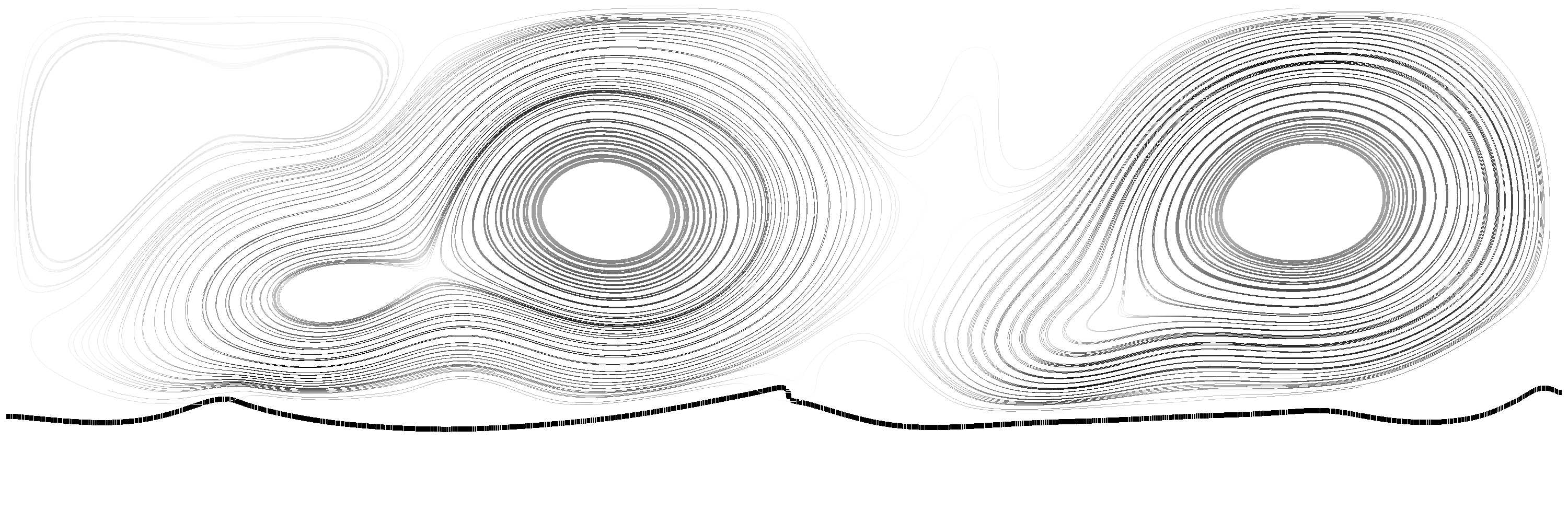}
  \hfill
  \\
  \hfill
  \includegraphics[width = 0.49\textwidth]{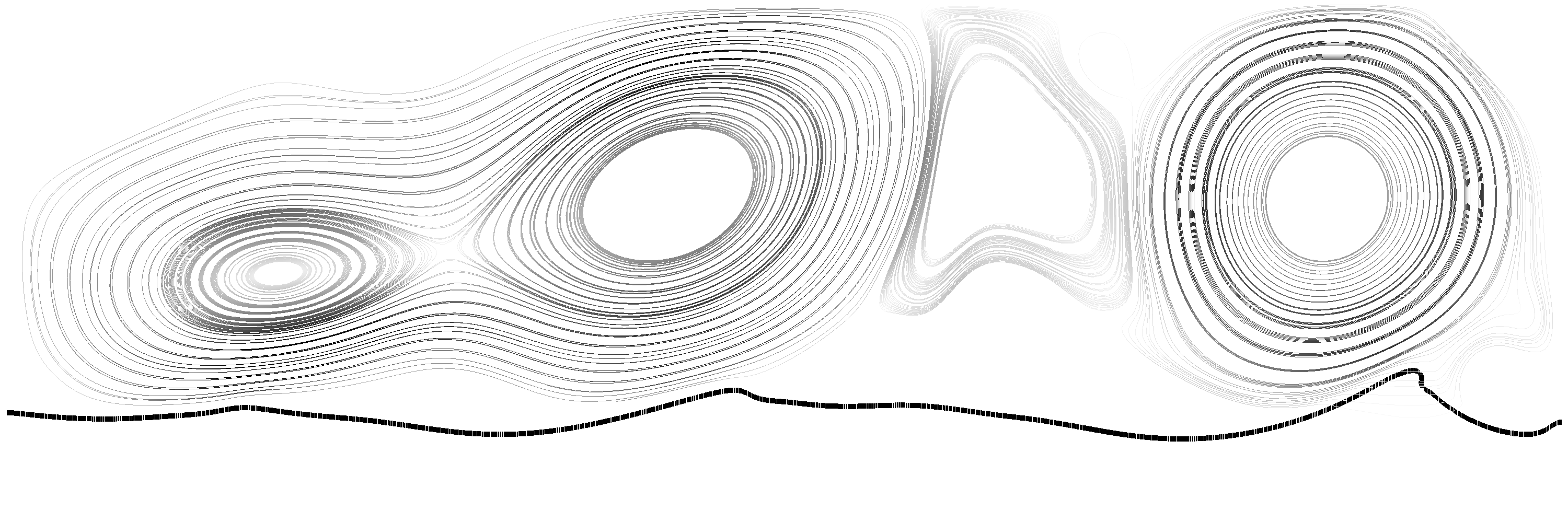}
  \hfill
  \includegraphics[width = 0.49\textwidth]{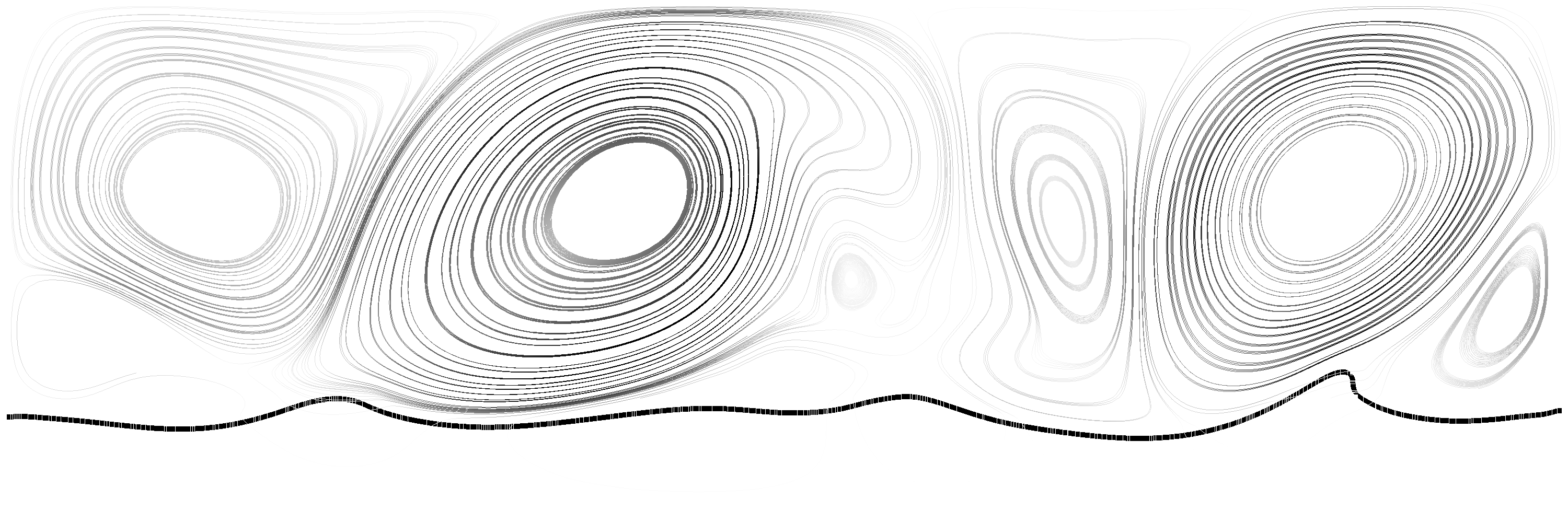}
  \hfill
  \\
  \hfill
  \includegraphics[width = 0.49\textwidth]{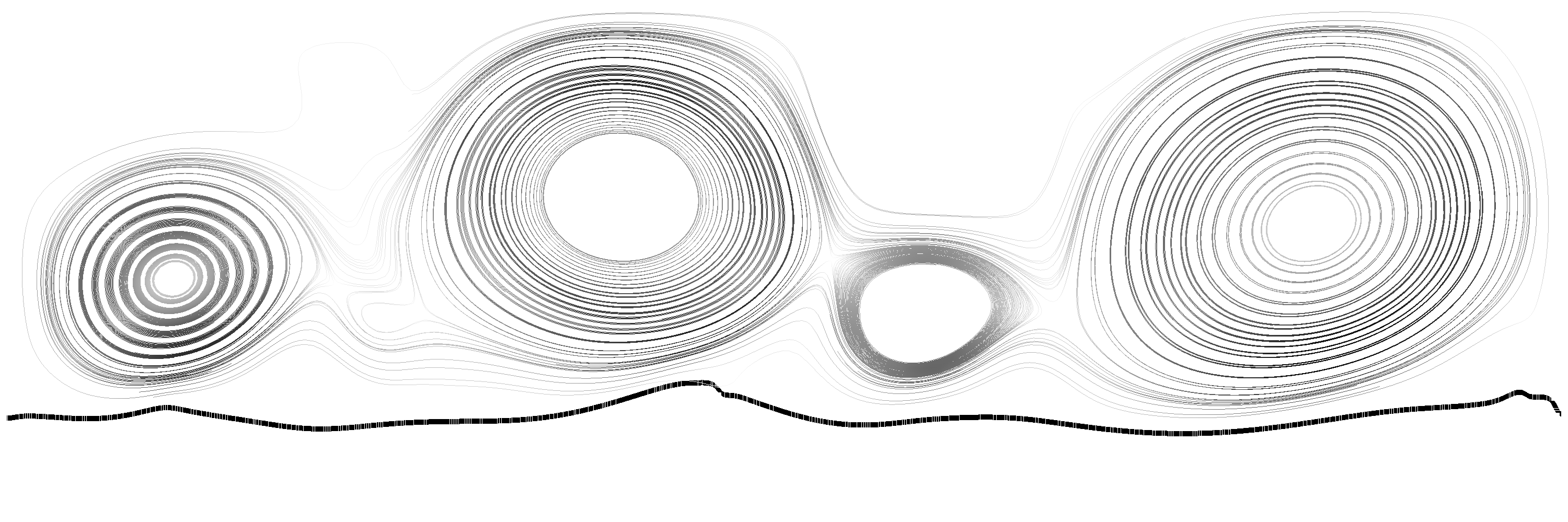}
  \hfill
  \includegraphics[width = 0.49\textwidth]{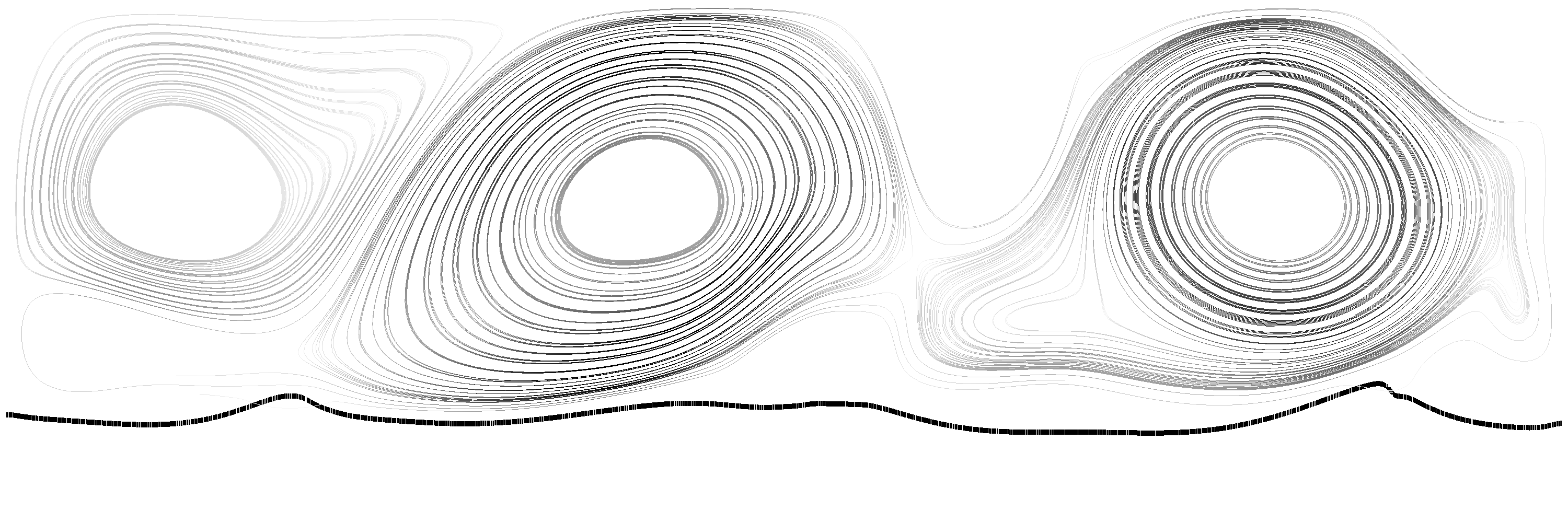}
  \hfill
  \\
  \caption{Snapshots of the evolution of $\varphi$ and $v$.
   For $t \in  \{1.7,3.3,5.0,6.7,8.3,10.0$\} 
   (left top to right bottom) 
   we present streamlines of $v$ in grayscale together with the zero level line
   of $\varphi$ in black. Darker streamlines means higher velocity. Due to the
   unphysical boundary data and the given forcing we observe large vortices
   that generate several waves at different locations. We stress, that
   especially breaking waves are captured by our approach as it is able to capture topological changes.
   }
  \label{m6:fig:num:results}
\end{figure}

In Fig. \ref{m6:fig:num:results} we show snapshots of the evolution of the
interface between water and air, given by the zero level line of $\varphi$, and
the velocity field  presented by streamlines of $v$, colored by $|v|$. We
observe that, despite the unphysical parameters and boundary data, the method
is able to deal with the complex two-phase interaction at the air-water
interface.

\section{Outlook on the direction of research}
\label{m6:sec:outlook}
The numerical results proposed in Sec. \ref{m6:sec:Numerics} are only preliminary and should be regarded as a proof of concept for the proposed diffuse interface approach. The method is able to cope with the
complex phenomena at the air-water interface. Further research is necessary to further develop our approach, and to make it applicable for real world scenarios. This includes

\bigskip
\noindent{\bf Boundary data}\\
As a first step the application of periodic boundary data for $v$ parallel to
the water surface will be incorporated together with an open boundary on the top and the
bottom of the domain. For $\varphi$ and $\mu$ periodic boundary conditions are sufficient in the water parallel directions only.

\bigskip
\noindent{\bf 3D computations}\\
For 3D computations an efficient solution of the linear systems arising
troughout the simulation is essential. Here results on preconditioning of the Cahn--Hilliard system from \cite{m6:BoschStollBenner_fastSolutionCH}  or a
multigrid approach as proposed in \cite{m6:KayWelford_multigridsolver} might be
used. For the solution of the Navier--Stokes equation well developed
preconditioners exist and we refer to
\cite{m6:Benzi_numericalSaddlePoint,m6:kayLoghinWelford_FpPreconditioner}.

\bigskip
\noindent{\bf Real world parameter}\\
Incorporating real world parameters will require several changes on the
architecture of the solver. Especially we note, that this will lead to large Reynolds number which require stabilization techniques like grad-div stabilization, which have to be encorporated into the finite element code. We note that the drawback of grad-div stabilization, namely a stronger coupling of the unknows,
does not appear here, as in \eqref{m6:eq:CHNSstrong1} all variables are coupled
anyway due to the term $2Dv = \nabla v + (\nabla v)^\top$.
 
\bigskip
\noindent{\bf Incorporation into Earth System Models}\\
If our concept proves applicable for the numerical simulation of the air-water region in atmosphere and ocean on the meter scale, it has to be incorporated into simulations on the next coarser (kilometer) scale. In this context homogenization concepts might be an option. We refer to e.g. \cite{m6:Eck2005}, where homogenization of phase field models has be done in a different context.    

\section*{Acknowledgement} The authors are greatful for many discussions with Jeff Carpenter from the Helmholtz Center in Geesthacht on practical issues related to the wind-wave coupling at the interface of atmosphere and ocean. The second author acknowledges support of the TRR 181 funded by the German Research Foundation.
%%%%%%%%%%%%%%%GKH-Literatur%%%%%%%%%%%%%%
\bibliographystyle{plain}
\bibliography{Quellen}

\end{document}